\documentclass[12pt,reqno]{amsart}

\usepackage{amsmath,amssymb,amsfonts,amsthm}
\numberwithin{equation}{section}
\newtheorem{theorem}{Theorem}[section]
\newtheorem{prop}[theorem]{Proposition}
\newtheorem{lemma}[theorem]{Lemma}

\usepackage{tikz}
\usetikzlibrary{matrix,arrows.meta,positioning, calc}
\usepackage[a4paper,
top=24mm,bottom=24mm,left=26mm,right=26mm
]{geometry}
\usepackage{graphicx}
\usepackage{booktabs}
\usepackage{mathtools}
\usepackage{makecell}

\usepackage[
]{cochineal}
\usepackage[cochineal,vvarbb,
cochf,cochrho]{newtxmath}
\usepackage[%
cal=zapfc, calscaled=1.15,
scr=boondox
]{mathalpha}

\usepackage[numbers,sort]{natbib}
\usepackage[
bookmarks=true,
bookmarksdepth=subsection,
bookmarksnumbered=true,
breaklinks=true,
colorlinks=true,
linkcolor=blue,
citecolor=magenta
]{hyperref}

\usepackage{listings}
\usepackage{xcolor}
\usepackage{tcolorbox}

\newcommand{\sageoutstyle}[1]{%
    \noindent\hspace*{8em}%
    {\ttfamily\small\itshape\color{black!70}#1}%
}
\lstdefinestyle{sageinput}{
    language=Python,
    backgroundcolor=\color{white},           %
    basicstyle=\ttfamily\small,
    keywordstyle=\color{blue!80!black}\bfseries,
    commentstyle=\color{green!50!black}\itshape,
    stringstyle=\color{purple},
    frame=leftline,                    %
    framerule=1.5pt,                   
    rulecolor=\color{gray!30},         
    xleftmargin=1.5em,
    framexleftmargin=0.5em,
    showstringspaces=false,
    breaklines=true,
    escapeinside={@}{@}                %
}

\DeclareMathOperator{\rank}{rank}
\DeclareMathOperator{\adj}{adj}
\DeclareMathOperator{\Aut}{Aut}
\DeclareMathOperator{\Tr}{Tr}

\DeclareMathOperator{\Hom}{Hom}

\DeclareMathOperator{\Sk}{Sk}

\DeclareMathDelimiter{\Norm}{\mathord}{largesymbols}{"3E}{largesymbols}{"3E}

\DeclareMathOperator{\Pois}{\mathscr{P}}
\DeclareMathOperator{\Sing}{Sing}

\allowdisplaybreaks
\begin{document}

\title[]{Character Variety of the Five-Punctured Sphere
  and
  the Determinantal Quintic Hypersurface}
\author[]{Kazuhiro Hikami}
\address{Faculty of Mathematics,
  Kyushu University,
  Fukuoka 819-0395, Japan.}

\email{
  \texttt{khikami@gmail.com}
}
\date{October 7, 2026.}

\begin{abstract}
  Using skein-theoretic methods,
  we study
  the~$\mathrm{SL}_2(\mathbb{C})$ character variety of the  five-punctured sphere
  with arbitrary conjugacy class at each of the five punctures.
  We show that
  a  Fricke--Klein--Vogt-type relation
  is realized  as a symmetric determinantal hypersurface.
  We also give  a realization of the character variety in terms of 
  cluster variables.
  We explicitly derive the Poisson structure of the simple closed
  curves on the surface, 
  and
  prove  that cluster  mutations induce 
  automorphisms of the character variety.
\end{abstract}
\subjclass[2000]{}
\maketitle

\section{Introduction}

Skein algebras on surface~$\Sigma_{g,n}$  have received renewed
interest
in their role 
as
quantum  character variety of the surface.
As seen from the fact that the
character variety has a connection with
the Coulomb branch,
the conformal block,
the Riemann--Hilbert
problem,
the Higgs
bundle,
$P=W$
conjecture, the Teichm{\"u}ller space,
wall crossing, and so on,
the skein algebra on the punctured sphere possesses a rich algebraic  structure.

Though the character variety of the four-punctured sphere has been well
understood since the 1800s,
the case of the $n$-punctured sphere~$\Sigma_{0,n}$
remains  far from  fully understood.
In~\cite{CalliMazzo18a}, the action of the mapping class group
of the five-punctured sphere~$\Sigma_{0,5}$
was
extensively studied in order to classify  finite orbits of the two-variable
Garnier system.
In~\cite{Simpson18a},
identified was
a relationship between a hypersurface and
the character variety $\mathfrak{X}_{\Sigma_{0,5}}$ with
a specific choice of
conjugacy classes at the five punctures.
It was constructed based on the pants decomposition of the punctured sphere.
See also~\cite{FassaLoray25a,Szabo23a,Girand16a} for related topics.

In this article,
motivated by~\cite{Simpson18a}, we study  the character variety of
the five-punctured sphere with arbitrary boundary conjugacy classes.
Using skein-theoretic methods together with cluster algebras,
we show
that
it
has a determinantal hypersurface expression.
See~\cite{Beauv00a,Piont06a} for  general properties  of determinantal
hypersurfaces
(see also~\cite{Dolga12Book}).

This paper is organized as follows.
In Section~\ref{sec:character}, we
study the character variety of the five-punctured sphere.
We obtain 
a symmetric determinantal hypersurface expression.
Theorem~\ref{thm:character} is our main result.
We then discuss the singular locus of the hypersurface.
In  support of Theorem~\ref{thm:character},
we analyze  the cluster algebraic structure of the simple closed curves
on~$\Sigma_{0,5}$
in Section~\ref{sec:cluster}.
We show that the  curves
can be written in terms of the cluster $y$-variables
following Bonahon--Wong~\cite{BonaWong16b}.
The Poisson algebra induced by the cluster variables reproduce  the
Goldman--Poisson algebra of the simple closed curves on~$\Sigma_{0,5}$.
We also show  that the cluster mutations induce  automorphisms 
of~$\Sigma_{0,5}$.
The last section is devoted to the concluding remarks.

\section{Character Variety and Determinantal Hypersurface}
\label{sec:character}
\subsection{Character Variety of $\Sigma_{0,5}$}

We define the $\mathrm{SL}_2(\mathbb{C})$ character variety
of a punctured surface~$\Sigma_{g,n}$
\begin{equation}
  \mathfrak{X}_{\Sigma_{g,n}}
  =
  \left.
    \Hom\left(
      \pi_1(\Sigma_{g,n}), \mathrm{SL}_2(\mathbb{C})
    \right)
    \middle/ \mathrm{SL}_2(\mathbb{C})
  \right.
  .
\end{equation}
For the $n$-punctured sphere,
the fundamental group $\pi_1(\Sigma_{0,n})$
is
\begin{equation}
  C_0 \,  C_1   \cdots  C_{n-1}=I   ,
\end{equation}
where $C_a$ denotes the monodromy of the $a$-th  boundary component.
It is known~\cite{RHorow72a}
that the character variety of the punctured sphere
is generated by
$\Tr (C_{i_1}C_{i_2} \cdots C_{i_\ell})$ for $ i_1< i_2 <\dots<
i_\ell $,
which is a consequence of
\begin{equation}
  \label{Trace_C}
  \Tr (C_1 \, C_2) + \Tr (C_1 \, C_2^{-1}) = \Tr(C_1) \, \Tr(C_2) ,
\end{equation}
for $C_i \in \mathrm{SL}_2(\mathbb{C})$.

This relation~\eqref{Trace_C} can be viewed as a classical skein relation.
The skein algebra $\Sk_A(\Sigma_{g,n})$ on $\Sigma_{g,n}$ is generated
by  framed links in the thickened surface $\Sigma_{g,n}\times [0,1]$
under the skein relation,
\begin{gather}
  \label{skein}
  \vcenter{\hbox{
    \begin{tikzpicture}%
      \draw [line width=1.2pt](0,1) --( 1,0) ;
      \draw[line width=10pt,white](0,0)--(1,1);
      \draw[line width=1.2pt](0,0)--(1,1);
    \end{tikzpicture}
  }}
  =
  A \,
  \vcenter{\hbox{
    \begin{tikzpicture}%
      \draw [line width=1.2pt] (1,1) to[out=-135,in=135] (1,0);
      \draw [line width=1.2pt] (0,1) to[out=-45,in=45] (0,0);
    \end{tikzpicture}
  }}
  +A^{-1} \,
  \vcenter{\hbox{
    \begin{tikzpicture}
      \draw [line width=1.2pt] (1,1) to[out=-135,in=-45] (0,1);
      \draw [line width=1.2pt] (1,0) to[out=135,in=45] (0,0);
    \end{tikzpicture}
  }}
  ,
  \\[2mm]
  \vcenter{\hbox{
    \begin{tikzpicture}%
      \draw [line width=1.2pt] (0,0) circle (0.5) ;
    \end{tikzpicture}
  }}
  = -A^2 - A^{-2} .
  \notag
\end{gather}
A multiplication~$\mathbb{x} \, \mathbb{y}$
of two links~$\mathbb{x}$ and~$\mathbb{y}$
is defined by placing~$\mathbb{x}$  vertically above~$\mathbb{y}$:
\begin{gather*}
  \mathbb{x} \, \mathbb{y}
  =
  \newcolumntype{C}{>{$}c<{$}}
  \begin{tabular}{|C|}
    \hline
    \hphantom{aa}\mathbb{x}\hphantom{aa}
    \\
    \hline
    \mathbb{y}
    \\
    \hline
  \end{tabular}
\end{gather*}
When $A=\pm 1$, the algebra is commutative,
and 
it is known~\cite{Bullo97a,PrzytSikor19a}
that $\Sk_{A=\pm 1}(\Sigma_{g,n})$ is isomorphic to 
the coordinate ring of the $\mathrm{SL}_2(\mathbb{C})$ character
variety~$\mathfrak{X}_{\Sigma_{g,n}}$
of the surface $\Sigma_{g,n}$.

\begin{figure}[htbp]
  \centering
  \begin{minipage}{0.3\linewidth}
    \centering
    \includegraphics[scale=0.9]{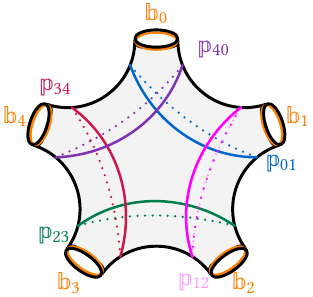}
  \end{minipage}
  \hfill
  \begin{minipage}{0.3\linewidth}
    \centering
    \includegraphics[scale=0.9]{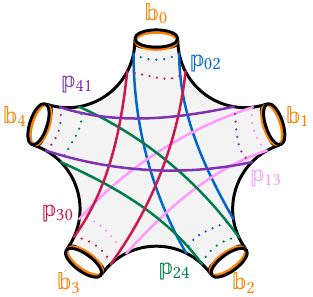}
  \end{minipage}
  \hfill
  \begin{minipage}{0.3\linewidth}
    \centering
    \includegraphics[scale=0.9]{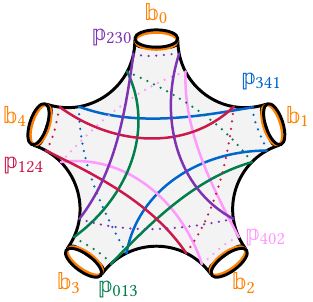}
  \end{minipage}
  \caption{Simple closed curves on $\Sigma_{0,5}$.}
  \label{fig:curves_Sigma05}
\end{figure}

On the five-punctured sphere
$\Sigma_{0,5}$,
we use  the simple closed curves shown
in Fig.~\ref{fig:curves_Sigma05}.
The curve~$\mathbb{b}_a$ denotes the boundary circle around the
$a$-th
puncture.
We denote   by $\mathbb{p}_{a,a+1}$   the curve encircling the
$a$-th and the $(a+1)$-th punctures where $a\in\{0,1,2,3,4\}$.
Throughout this article,  the subscript
should be read modulo $5$.
We also use  curves $\mathbb{p}_{a,a+2}$
(resp. $\mathbb{p}_{a,a+1,a+3}$) as the curves enclosing
the        punctures $a$ and $a+2$
(resp. $a$, $a+1$ and $a+3$) as in Fig.~\ref{fig:curves_Sigma05}.
Hereafter
we denote by
$p_\bullet$ (resp. $b_a$)
the element of the classical commutative algebra $\Sk_{A=1}(\Sigma_{0,5})$
corresponding to the curve~$\mathbb{p}_\bullet $
(resp. $\mathbb{b}_a $).
For brevity,
we write
a generating set  of $\Sk_{A=1}(\Sigma_{0,5})$
as
\begin{equation}
  \label{20-curve}
  \mathbf{c}
  =
  \begin{bmatrix}
    p_{01} & p_{13} & p_{402} & b_0 \\
    p_{12} & p_{24} & p_{013} & b_1 \\
    p_{23} & p_{30} & p_{124} & b_2 \\
    p_{34} & p_{41} & p_{230} & b_3 \\
    p_{40} & p_{02} & p_{341} & b_4     
  \end{bmatrix}
  .
\end{equation}

\begin{figure}[htbp]
  \centering
  \includegraphics[scale=.9]{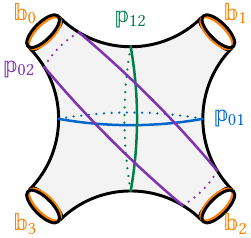}  
  \caption{A set of simple closed curves on $\Sigma_{0,4}$ is denoted.
  The boundary  circles are $\mathbb{b}_a$.}
  \label{fig:curves_Sigma04}
\end{figure}

The character variety $\mathfrak{X}_{\Sigma_{0,4}}$ of
the four-punctured sphere is intimately related to  Painlev{\'e}~VI
and  has been extensively studied.
In particular,
for the simple closed curves in Fig.~\ref{fig:curves_Sigma04},
we have the
Fricke--Klein--Vogt
relation~\cite{Vogt89a,FrickeKleinBook1897,Iwasaki03a,Goldman09a}
\begin{equation}
  \label{Fricke-Klein=0}
  W_{{b}_0,{b}_1,{b}_2,{b}_3}
  ({p}_{01}, {p}_{12}, {p}_{20})
  =0 ,
\end{equation}
where
\begin{multline}
  \label{Fricke-Klein}
  W_{{b}_0,{b}_1,{b}_2,{b}_3}
  ({p}_{01}, {p}_{12}, {p}_{20})
  =
  -  {p}_{01}  {p}_{12}  {p}_{20}
  +
  {p}_{01}^2 + {p}_{12}^2 + {p}_{20}^2
  +
  \left({b}_0 {b}_1 + {b}_2 {b}_3\right)
  {p}_{01}
  \\
  +
  \left({b}_1 {b}_2 + {b}_3 {b}_0\right)
  {p}_{12}
  +
  \left({b}_2 {b}_0 + {b}_1 {b}_3\right)
  {p}_{20}
  +
  {b}_0^2+  {b}_1^2+  {b}_2^2+  {b}_3^2
  +  {b}_0{b}_1{b}_2{b}_3
  -4 .
\end{multline}

The character variety~$\mathfrak{X}_{\Sigma_{0,5}}$ and the skein
algebra~$\Sk_A(\Sigma_{0,5})$
have a more complicated structure.
See~\cite{CookeLacab22a} for the presentation
of~$\Sk_A(\Sigma_{0,5})$.
To have a nice presentation at the classical limit,
we recall that
several subsurfaces $\Sigma_{0,4}$ of $\Sigma_{0,5}$
give rise to the Fricke--Klein--Vogt relations as
(see also~\cite{CalliMazzo18a})
\begin{equation}
  \label{W_for_2-p}
  W_{{b}_a,{b}_{a+1},{b}_{a+2},{p}_{a+3,a+4}}
  ({p}_{a,a+1}, {p}_{a+1,a+2}, x)=0,
  \quad \text{for $x={p}_{a,a+2}$,
    ${p}_{a+3,a+4,a+1}$,}
\end{equation}
which gives
\begin{equation}
  \label{pp-1}
  {p}_{a,a+2}+{p}_{a+3,a+4,a+1}
  =
  {p}_{a,a+1}\, {p}_{a+1,a+2}
  -
  \left(
    {b}_a \, {b}_{a+2} + {b}_{a+1} \,    {p}_{a+3,a+4}
  \right) .
\end{equation}

Additional gluing relations  are required.
The skein relation yields
\begin{align}
  \label{eq:6}
  \mathbb{p}_{02} \, \mathbb{p}_{23}
  &=
  \raisebox{-0.5\height}{
    \includegraphics[scale=0.5]{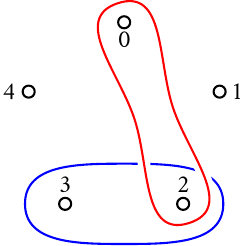}
  }
  =
  A^2 \,
  \raisebox{-0.5\height}{
    \includegraphics[scale=0.5]{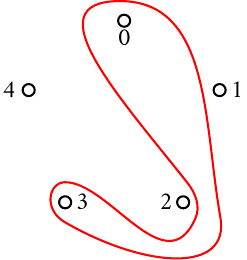}
  }
  +A^{-2} \, \mathbb{p}_{30} + \mathbb{b}_0\, \mathbb{b}_{3}
  +
  \mathbb{b}_2\, \mathbb{p}_{230} ,
  \\
  \mathbb{p}_{402} \, \mathbb{p}_{01}
  &=
  \raisebox{-0.5\height}{
    \includegraphics[scale=0.5]{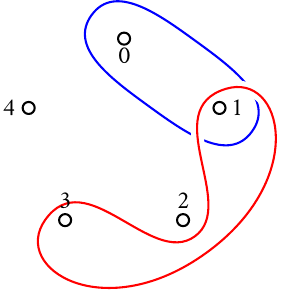}
  }
  =
  A^2 \,
  \raisebox{-0.5\height}{
    \includegraphics[scale=0.5]{pentagon_skein3}
  }
  +A^{-2} \, \mathbb{p}_{124} + \mathbb{b}_0\, \mathbb{b}_{3}
  +
    \mathbb{b}_1\, \mathbb{p}_{24} ,
\end{align}
where
white circles denote  punctures on the sphere.
By subtracting the two identities and taking
the classical limit~$A\to 1$,
we obtain
\begin{equation}
  \label{pp-2}
  {p}_{402} \, {p}_{01}
  +{p}_{30}+{b}_2\, {p}_{230}
  =
  {p}_{02} \, {p}_{23}
  +{p}_{124}+{b}_1\, {p}_{24} .
\end{equation}
Using~\eqref{pp-1},
we eliminate~${p}_{a,a+1,a+3}$ from~\eqref{pp-2}.
By cyclic rotation,
we obtain
five equations 
for ${p}_{a,a+2}$ as
\begin{equation}
  \label{our-M-eq}
  \mathbf{M}
  \begin{pmatrix}
    {p}_{13} \\
    {p}_{24} \\
    {p}_{30}\\
    {p}_{41} \\
    {p}_{02} 
  \end{pmatrix}
  =
  \begin{pmatrix}
    -{p}_{34} \, {p}_{40} \, {p}_{01}    + {p}_{23} \, {p}_{12} +
    {b}_4 \, {b}_1  \, {p}_{34}    + {b}_3 \, {b}_0 \,{p}_{01}
    - {b}_2 \, {p}_{40} +    {b}_2 \, {b}_4  \, {b}_0    -{b}_3 \, {b}_1
    \\
    -{p}_{40} \, {p}_{01} \, {p}_{12}    + {p}_{34} \, {p}_{23} +
    {b}_0 \, {b}_2  \, {p}_{40}    + {b}_4 \, {b}_1 \,{p}_{12}
    - {b}_3 \, {p}_{01} +    {b}_3 \, {b}_0  \, {b}_1
    -{b}_4 \, {b}_2
    \\
    -{p}_{01} \, {p}_{12} \, {p}_{23}    + {p}_{40} \, {p}_{34} +
    {b}_1 \, {b}_3  \, {p}_{01}    + {b}_0 \, {b}_2 \,{p}_{23}
    - {b}_4 \, {p}_{12} +    {b}_4 \, {b}_1  \, {b}_2
    -{b}_0 \, {b}_3
    \\
        -{p}_{12} \, {p}_{23} \, {p}_{34}    + {p}_{01} \, {p}_{40} +
    {b}_2 \, {b}_4  \, {p}_{12}    + {b}_1 \, {b}_3 \,{p}_{34}
    - {b}_0 \, {p}_{23} +    {b}_0 \, {b}_2  \, {b}_3
    -{b}_1 \, {b}_4
    \\
    -{p}_{23} \, {p}_{34} \, {p}_{40}    + {p}_{12} \, {p}_{01} +
    {b}_3 \, {b}_0  \, {p}_{23}    + {b}_2 \, {b}_4 \,{p}_{40}
    - {b}_1 \, {p}_{34} +    {b}_1 \, {b}_3  \, {b}_4
    -{b}_2 \, {b}_0
  \end{pmatrix} ,
\end{equation}
where the symmetric matrix $\mathbf{M}$ is given by
\begin{align}
  \mathbf{M}
  & =
    \mathbf{M}_{{b}_0,{b}_1,{b}_2,{b}_3,{b}_4}
    ({p}_{01},
    {p}_{12},{p}_{23},{p}_{34},{p}_{40}
    )
    \nonumber
  \\
  &=
    \begin{pmatrix}
      2 & - {b}_0 & - {p}_{01}
      & - {p}_{34}  & - {b}_4
      \\
      - {b}_0 & 2 & - {b}_1
      & - {p}_{12}    & - {p}_{40}
      \\
      - {p}_{01} & - {b}_1 & 2
      & - {b}_2   & - {p}_{23}
      \\
      - {p}_{34} & - {p}_{12} & - {b}_2
      & 2   & - {b}_3
      \\
      - {b}_4 & - {p}_{40} & - {p}_{23}
      &- {b}_3 & 2      
    \end{pmatrix} .
    \label{define_Q}
\end{align}
If $\mathbf{M}$ were full rank,~\eqref{our-M-eq} would admit a unique solution,
and~\eqref{pp-1} would then force
$
{p}_{a,a+2}={p}_{a+3,a+4,a+1}
$ for all $a$,
which is not true in general.
We therefore  conclude that
$\mathbf{p}=(p_{01}, p_{12}, p_{23}, p_{34}, p_{40})$
has a determinantal  hypersurface expression,
\begin{equation}
  \label{det-Q}
  \det \mathbf{M}=0 .
\end{equation}
Explicitly,
this is a quintic hypersurface:
\begin{multline}
  \label{det-Q-explicit}
  - {p}_{01} \, {p}_{12} \, {p}_{23} \,
  {p}_{34} \, {p}_{40}
  +\sum_a
  \left(
    {p}_{a,a+1} \, {p}_{a+1,a+2}
    -{b}_{a}  \, {b}_{a+2}
    - {b}_{a+1} \, {p}_{a+3,a+4}
  \right)^2
  \\
  -4 \sum_a {p}_{a,a+1}
  \left(
    {p}_{a,a+1}
    +{b}_{a+2}      {p}_{a+3,a+4} 
    + {b}_a {b}_{a+1} +
    {b}_{a+2}       {b}_{a+3} {b}_{a+4}
  \right) 
  \\
  +\sum_a {b}_{a+4}
  \left(
    {p}_{a,a+1} \, {p}_{a+1,a+2}
    - {b}_{a+1} \, {p}_{a+3,a+4}
    -{b}_{a}  \, {b}_{a+2}
  \right){p}_{a+1,a+2} {p}_{a+2,a+3}
    +\sum_a
  {b}_a^2 {b}_{a+2} {b}_{a+3} {p}_{a+2,a+3}
  \\
  - 4 \sum_a {b}_a^2
  - {b}_0 \, {b}_1 \, {b}_2 \, {b}_3 \,
  {b}_4 +16
  =0 .
\end{multline}
One sees that,
for~$\mathbf{b}=(0,0,0,0,0)$,~\eqref{det-Q-explicit}
reduces 
to
\begin{multline*}
  - {p}_{01} \, {p}_{12} \, {p}_{23} \,
  {p}_{34} \, {p}_{40}
  + \left( {p}_{01} \, {p}_{12}\right)^2
  + \left( {p}_{12} \, {p}_{23}\right)^2
  + \left( {p}_{23} \, {p}_{34}\right)^2
  + \left( {p}_{34} \, {p}_{40}\right)^2
  + \left( {p}_{40} \, {p}_{01}\right)^2
  \\
  -4 \left( {p}_{01}^2 + {p}_{12}^2
    + {p}_{23}^2 + {p}_{34}^2 + {p}_{40}^2
  \right)
  +16
  =0 ,
\end{multline*}
which agrees with the result in~\cite{Simpson18a}.
Therein used was the Fenchel--Nielsen coordinates to a pants
decomposition of~$\Sigma_{0,5}$.

We should note that,
using the cluster variables~\eqref{boundary_cluster}
and~\eqref{p_a_a+1_cluster} which will be discussed later,
we can check  by \texttt{Mathematica}
\begin{equation}
  \label{rank-Q}
  \rank \mathbf{M}=4 ,
\end{equation}
at general points of the moduli space.

\begin{theorem}
  \label{thm:character}
  The character variety $\mathfrak{X}_{\Sigma_{0,5}}$ of
  the five-punctured sphere
  is naturally embedded in
  $\mathbb{A}^{20}$ with coordinates~$\mathbf{c}$~\eqref{20-curve}
  satisfying~\eqref{W_for_2-p},~\eqref{our-M-eq}
  with the symmetric rank~$4$ matrix $\mathbf{M}$~\eqref{define_Q}.  
\end{theorem}

The skein relation~\eqref{skein} on $\Sigma_{0,5}$
also proves the following.

\begin{prop}
  The kernel of~$\mathbf{M}$~\eqref{define_Q} is spanned by
  $\mathbf{v}  =
  (v_0, \cdots, v_4)^\top $
  defined by
\begin{equation}
  \label{vector_kernel}
  \mathbf{v}
  =
  \begin{pmatrix}
    {p}_{13}-{p}_{402} \\
    {p}_{24}-{p}_{013} \\
    {p}_{30}-{p}_{124} \\
    {p}_{41}-{p}_{230} \\
    {p}_{02}-{p}_{341} 
  \end{pmatrix}
  \in \ker \mathbf{M} .
\end{equation}
\end{prop}

We note that~\eqref{vector_kernel}
can also be checked by \texttt{Mathematica} using the cluster
variables~\eqref{boundary_cluster},~\eqref{p_a_a+1_cluster},
and~\eqref{p_a_a+2_cluster}.
Since~\eqref{rank-Q} implies
$\rank (\adj(\mathbf{M}))=1$,
the symmetric adjugate matrix $\adj(\mathbf{M})$ is proportional to the outer
product of~$\mathbf{v}$~\eqref{vector_kernel} with itself.
In fact using the same cluster
coordinates,
we  find that
\begin{equation}
  \label{adj-Q_v_v}
  \adj(\mathbf{M})
  =
  \mathbf{v} \, \mathbf{v}^\top .
\end{equation}

\subsection{Determinantal Quintic Hypersurface}

We denote the determinantal quintic hypersurface
in the coordinates~$\mathbf{p}=(p_{01}, p_{12}, p_{23}, p_{34}, p_{40})\in
\mathbb{A}^5$
by
\begin{equation}
\mathcal{Q}=
\left\{
  \mathbf{p} \in \mathbb{A}^5
  ~\middle|~
  \det \mathbf{M}(\mathbf{p})=0
\right\} .
\end{equation}
and set~\cite{Simpson18a}
\begin{equation}
  \phi:
  \mathfrak{X}_{\Sigma_{0,5}}
  \to \mathcal{Q} .
\end{equation}
The singular locus of~$\mathcal{Q}$
is
given by
\begin{equation}
  \label{Q_critical}
  \frac{\partial}{\partial {p}_{a,a+1}} \det \mathbf{M}
  =\Tr \left( \adj(\mathbf{M}) \,
    \frac{\partial \mathbf{M}}{\partial {p}_{a,a+1}}
  \right)
  =0 ,
\end{equation}
which yields
$v_a \, v_{a+2}=0$
for all~$a$
due to~\eqref{adj-Q_v_v}.
Thus we have $\bigcup_a \{ v_a=v_{a+1}=v_{a+2}=0 \}$.
When $v_a=v_{a+1}=v_{a+2}=0$,~\eqref{vector_kernel}
gives
\begin{equation*}
\begin{pmatrix}
  2 & -b_{a+3} \\
  -b_{a+3} & 2
\end{pmatrix}
\begin{pmatrix}
  v_{a+3} \\
  v_{a+4}
\end{pmatrix}
=
\mathbf{0} ,
\end{equation*}
and the singular locus is 
\begin{equation}
  \label{singular_locus_b}
  \{ \mathbf{v}=\mathbf{0} \}
  \cup \bigcup_a
  \left\{
    {b}_a=\pm 2,
    v_a=\pm v_{a+1} ,
    v_{a+2}=v_{a+3}=v_{a+4}=0  
  \right\} .
\end{equation}
The case $\mathbf{v}=\mathbf{0}$ corresponds to the essential
singularity whose corank is greater than or equal to~$2$.
The accidental singularities occur when
$b_a=\pm 2$ for some $a$.

At a generic point,
the singular locus
of $\mathcal{Q}$ is
2-dimensional, and we expect that
\begin{equation}
  \label{Singular_rank_3}
  \Sing( \mathcal{Q})=
  \left\{
    \mathbf{p} \in \mathbb{A}^5
    ~\middle|~
    \rank \mathbf{M}(\mathbf{p})=3
  \right\} .
\end{equation}
As
$\mathbf{M} \sim
\begin{psmallmatrix}
  I_3 & O \\
  O & O
\end{psmallmatrix}
$,
the singularity is  of type~$A_1$ with
exceptional fiber~$\mathbb{P}^1$.

The Fricke--Klein--Vogt relation~\eqref{W_for_2-p}
proves that
the condition $v_{a-1}=0$, \emph{i.e.},
${p}_{a,a+2}={p}_{a+3,a+4,a+1}$,
means that
\begin{multline}
  \label{eq:22}
  \left( {b}_a \, {b}_{a+2} +
    {b}_{a+1} \, {p}_{a+3,a+4}
    - {p}_{a,a+1} \,    {p}_{a+1,a+2}
  \right)^2
  \\
  =
  4 \Bigl\{
  {p}_{a,a+1}^2 + {p}_{a+1,a+2}^2 +
  {p}_{a+3,a+4}^2 +
  \left(
    {b}_a \, {b}_{a+1}
    + {b}_{a+2} \, {p}_{a+3,a+4}
  \right) {p}_{a,a+1}
  +
  \left(
    {b}_{a+1}\, {b}_{a+2}
    +{b}_a \,    {p}_{a+3,a+4}
  \right)  {p}_{a+1,a+2}
  \\
  + {b}_a \, {b}_{a+1} \, {b}_{a+2}  \, {p}_{a+3,a+4}
  +
  {b}_a^2 + {b}_{a+1}^2 +{b}_{a+2}^2
  -4
  \Bigr\} .
\end{multline}
The coincidence between the ramification locus of~$\phi$ and the
singular locus of $\mathcal{Q}$ was also  noticed in~\cite{Simpson18a}.

Let~$\widebar{\mathcal{Q}}\subset \mathbb{P}^5$ denote
the  projective hypersurface
given by the homogenization of
$\mathbf{M}$~\eqref{define_Q} with
${p}_{a,a+1}=P_{a,a+1}/P_{00}$.
At infinity~$P_{00}=0$, we have
$P_{01} P_{12} P_{23} P_{34} P_{40}=0$,
and
the divisor at infinity of $\widebar{\mathcal{Q}}$
is~\cite{Simpson18a}
\begin{equation}
  \widebar{\mathcal{Q}}_\infty
  =
  \bigcup_a \{
  P_{00}=
  P_{a,a+1}=0\} .
\end{equation}
See~\cite{Simpson16a,TehraFrohm26a,Komyo15a} for studies on the boundary
divisors of compactification of the character variety.
By the same computation with~\eqref{Q_critical}, 
we have 
\begin{multline}
  \label{singular_infinity}
  \Sing(\widebar{\mathcal{Q}}_\infty)
  =
  \bigcup_a
  \left\{
    P_{00}=
    P_{a,a+1}=P_{a+1,a+2}=P_{a+3,a+4}=0
  \right\}
  \\
  \cup
  \left\{
    P_{00}=P_{a-1,a}=P_{a,a+1}=
    P_{a+1,a+2}^2+ b_{a}  \, P_{a+1,a+2} \, P_{a+3,a+4}+P_{a+3,a+4}^2=0
  \right\} ,
\end{multline}
which degenerates when $b_{a}=\pm 2$.
It should be  noted that~\eqref{singular_infinity}
was studied  in~\cite{Simpson18a} for the case of $b_a=0$.

In the following,
we  examine specific  choices of~$\mathbf{b}$ separately.
We show that
the accidental singular loci in~\eqref{singular_locus_b}
give rise to  singularities of~$\widetilde{\mathcal{Q}}$, defined by
\begin{equation}
\pi: \widetilde{\mathcal{Q}} \to \mathcal{Q},  
\end{equation}
where~$\widetilde{\mathcal{Q}}$ is  the kernel  incidence variety
defined by
\begin{equation}
  \label{define_Q_tilde}
  \widetilde{\mathcal{Q}}
  =
  \left\{
    (
    \mathbf{p}, [\mathbf{w}])
    \in \mathbb{A}^5 \times \mathbb{P}^4
    ~\middle|~
    {\mathbf{M}}(\mathbf{p}) \, \mathbf{w}=
    \mathbf{0}
  \right\} .
\end{equation}
The singularity of $\widetilde{\mathcal{Q}}$
is determined by the
Jacobian~$\mathbf{J}$ of $\mathbf{M} \, \mathbf{w}$
with respect to the variables~$\mathbf{w}$ and~$\mathbf{p}$,
where $\mathbf{w}= \left(w_0,w_1,w_2,w_3,w_4 \right)^\top $:
\begin{equation}
  \label{Jacobian_Q-w}
  \mathbf{J}
  =
  \begin{pmatrix}
    2 & - {b}_0 & - {p}_{01}
    & - {p}_{34}  & - {b}_4
    & -w_2 & 0 & 0 & -w_3 & 0
    \\
    - {b}_0 & 2 & - {b}_1
    & - {p}_{12}    & - {p}_{40}
    & 0 & -w_3 & 0 & 0 & -w_4
    \\
    - {p}_{01} & - {b}_1 & 2
    & - {b}_2   & - {p}_{23}
    & -w_0 & 0 & -w_4 & 0 & 0
    \\
    - {p}_{34} & - {p}_{12} & - {b}_2
    & 2   & - {b}_3
    & 0 & -w_1 & 0 & -w_0 & 0
    \\
    - {b}_4 & - {p}_{40} & - {p}_{23}
    &- {b}_3 & 2
    & 0 & 0 & -w_2 & 0 & -w_1
  \end{pmatrix} .
\end{equation}
We need to check the condition for $\rank \mathbf{J}<5$.
Hereafter we denote by $\mathcal{S}_r$ the stratum
\begin{equation}
  \label{eq:16}
  \mathcal{S}_r=
  \Sing(\mathcal{Q}) \cap \{ \mathbf{p} ~|~ \rank \mathbf{M}(\mathbf{p})=r\} .
\end{equation}

\subsubsection{$\mathbf{b}=(0,0,0,0,0)$}
This is the
generic point studied in~\cite{Simpson18a}, and the entire singular locus is
\begin{equation*}
  \Sing(\mathcal{Q})=
  \mathcal{S}_3 ,
\end{equation*}
and
we have $\pi^{-1}(\mathbf{p})=\mathbb{P}^1$ for $\mathbf{p}\in \mathcal{S}_3$.

\subsubsection{$\mathbf{b}=(-2,0,0,0,0)$}

We have  an accidental singularity as seen from~\eqref{singular_locus_b},
and
\begin{equation*}
  \Sing(\mathcal{Q})=
  \mathcal{S}_4 \cup \mathcal{S}_3 \cup \mathcal{S}_2 ,
\end{equation*}
where
\begin{align*}
  \mathcal{S}_2
  & =
  \{ \pm (0,2,2,2,0),
    \pm (0, 2,-2, 2,0)\} .
\end{align*}
This set of  rank~$2$ singularity is induced from the accidental singular
locus,
\begin{equation*}
    \mathcal{S}_4
    =
    \left\{
      (0,s,t,s,0)
      ~\middle|~ (s^2-4)(t^2-4)\neq 0
    \right\} .
\end{equation*}
The Jacobian criterion for~\eqref{Jacobian_Q-w} shows that
the accidental singularity induces
a $\mathbb{P}^0$ fiber
\begin{equation}
  \Sing(\widetilde{\mathcal{Q}})
  =
  \left\{
    \left(
      (0,s, t, s, 0),
      [1:-1:0:0:0]
    \right)
  \right\} .
\end{equation}

\subsubsection{$\mathbf{b}=(-2,-2,0,0,0)$}
We have
\begin{equation*}
  \Sing(\mathcal{Q})=
  \mathcal{S}_4 \cup \mathcal{S}_3 ,
\end{equation*}
where the
two accidental singularities give
\begin{equation*}
  \mathcal{S}_4
  =
  \left\{ (-2,s,t,s,0)
    ~\middle|~
    s \, t \neq 0
  \right\}
  \cup
  \left\{ (-2,0,s,t,s)
    ~\middle|~
    4t^2 \neq s^2(t^2-4)
  \right\} .
\end{equation*}
The Jacobian~\eqref{Jacobian_Q-w}
indicates that
\begin{equation}
  \Sing(\widetilde{\mathcal{Q}})
  =
  \left\{
    \left(
      (-2,s,t,s,0), [1:-1:0:0:0]
    \right)
  \right\}
  \cup
  \left\{
    \left(
      (-2,0,s,t,s), [0:1:-1:0:0]
    \right)
  \right\},
\end{equation}
and that
the fiber jumps at the intersection of the two loci as
\begin{equation}
  \left. \Sing(\widetilde{\mathcal{Q}})
  \right|_{\mathbf{p}=(-2,0,0,0,0)}
  =
  \left\{
    w_0+w_1+w_2=w_3=w_4=0
  \right\} .
\end{equation}

\subsubsection{$\mathbf{b}=(-2,-2,-2,0,0)$}
We have
\begin{equation*}
  \Sing(\mathcal{Q})
  =
  \mathcal{S}_4 \cup  \mathcal{S}_3 \cup  \mathcal{S}_2,
\end{equation*}
where
\begin{equation*}
  \mathcal{S}_2= \{ (-2,-2,0,-2,0) \}.
\end{equation*}
The three accidental singular loci induce
\begin{multline}
  \Sing(\widetilde{\mathcal{Q}})
  =
  \left\{
    ( (-2,s,t,s,0), [1:-1:0:0:0])
  \right\}
  \\
  \cup
  \left\{
    ( (-2,-2,s,t,s), [0:1:-1:0:0])
  \right\} \cup
  \left\{
    ( (s,-2,0,s,t), [0:0:1:-1:0])
  \right\} .
\end{multline}
The three singular loci intersect at a point
$\mathbf{p}=(-2,-2,0,-2,0)$, where
the fiber jumps  to
\begin{multline}
  \left.
    \Sing(\widetilde{\mathcal{Q}})
  \right|_{\mathbf{p}=(-2,-2,0,-2,0)}
  =
  \left\{
    w_0=w_1+w_2+w_3=w_4=0
  \right\}
  \\
  \cup
  \left\{
    w_0+w_1=w_2+w_3=w_4=0
  \right\}
  \cup
  \left\{
    w_0+w_1+w_2=w_3=w_4=0
  \right\} .
\end{multline}

\subsubsection{$\mathbf{b}=(-2,-2,-2,-2,0)$}
We have
\begin{equation*}
  \Sing(\mathcal{Q})
  = \mathcal{S}_4 \cup \mathcal{S}_3,
\end{equation*}
and the four accidental singularities induce
\begin{multline}
  \Sing( \widetilde{\mathcal{Q}} )
  =
  \left\{
    ( (-2,s,t,s,0) , [1:-1:0:0:0])
  \right\}
  \cup
    \left\{
    ( (-2,-2,s,t,s) , [0:1:-1:0:0])
  \right\}
  \\
  \cup
    \left\{
    ( (s,-2,-2,s,t) , [0:0:1:-1:0])
  \right\}
  \cup
    \left\{
    ( (t,s,-2,0,s) , [0:0:0:1:-1])
  \right\} .
\end{multline}
The four singular loci intersect each other at distinct points,
and
the jumping fibers at each intersection
are as follows:
\begin{equation*}
  \begin{array}[]{cc}
    \toprule
    \mathbf{p}
    & \left. \Sing(\widetilde{\mathcal{Q}})
      \right|_{\mathbf{p}}
      \\
    \midrule
    (-2,-2,0,-2,0)
    &
    w_0+w_1+w_2=w_3=w_4=0
  \\
    (-2,-2,-2,-2,0)
    &
    w_0+w_1=w_2+w_3=w_4=0
  \\
    (-2,0,-2,0,0)
    &
    w_0+w_1=w_2=w_3+w_4=0
  \\
    (-2,-2,-2,-2,-2)
    &
    w_0=w_1+w_2+w_3=w_4=0
  \\
    (-2,-2,-2,0,-2)
    &
    w_0=w_1+w_2=w_3+w_4=0
  \\
    (0,-2,-2,0,-2)
  &
    w_0=w_1=w_2+w_3+w_4=0
    \\
    \bottomrule
  \end{array}
\end{equation*}

\subsubsection{$\mathbf{b}=(-2,-2,-2,-2,-2)$}

As in the case of $\Sigma_{0,4}$ where
the Fricke--Klein--Vogt relation~\eqref{W_for_2-p} is most singular at
$b_a=\pm 2$ for~$\forall a$,
this is the most singular
case of~$\Sigma_{0,5}$.
We have a rank~$1$ singularity
\begin{equation*}
  \Sing (\mathcal{Q})=
  \mathcal{S}_4 \cup \mathcal{S}_3 \cup \mathcal{S}_1 ,
\end{equation*}
where
\begin{equation*}
  \mathcal{S}_1
  =
  \{(-2,-2,-2,-2,-2)\}.  
\end{equation*}
At $\mathbf{p}=(-2,-2,-2,-2,-2)$
the multiplicity is 4, since,
when we write
$p_{a,a+1}=-2+x_{a}$,
the hypersurface is given by
\begin{equation*}
  - x_{0} \,  x_{1} \,  x_{2} \,  x_{3} \,  x_{4}
  +\sum_a  \left(
    x_{a}^2 \, x_{a+2}^2
    -2 \, x_{a-1} \, x_{a}^2 \, x_{a+1}
    +2 \, x_{a} \, x_{a+1} \, x_{a+2} \,x_{a+3}
  \right)
  =0 .
\end{equation*}

\begin{figure}[htbp]
  \centering
    \includegraphics[scale=.9]{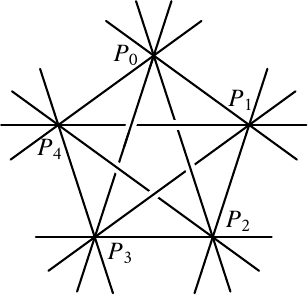}
    \caption{$E_{a,b}$~\eqref{line_L_a_b} in $\mathbb{P}^4$ pass through 
      $P_a$ and $P_b$.
      Here $P_0=[1:-1:0:0:0]$, $P_1=[0:1:-1:0:0]$,
      $P_2=[0:0:1:-1:0]$, $P_3=[0:0:0:1:-1]$, and
      $P_4=[-1:0:0:0:1]$.
    }
  \label{fig:pentagram}
\end{figure}

We have five singular loci as
\begin{multline}
  \Sing( \widetilde{\mathcal{Q}} )
  =
  \left\{
    ( (-2,s,t,s,-2) , [1:-1:0:0:0])
  \right\}
  \cup
  \left\{
    ( (-2,-2,s,t,s) , [0:1:-1:0:0])
  \right\}
  \\
  \cup
  \left\{
    ( (s,-2,-2,s,t) , [0:0:1:-1:0])
  \right\}
  \cup
  \left\{
    ( (t,s,-2,-2,s) , [0:0:0:1:-1])
  \right\}
  \\
  \cup
  \left\{
    ( (s,t,s,-2,-2) , [-1:0:0:0:1])
  \right\} .
\end{multline}
They intersect at a single point $\mathbf{p}=(-2,-2,-2,-2,-2)$,
and
the fiber jumps to
\begin{equation}
  \label{eq:9}
  \left.
    \Sing(\widetilde{\mathcal{Q}})
  \right|_{\mathbf{p}=(-2,-2,-2,-2,-2)}
  =
  \bigcup_{a} E_{a,a+1}\cup E_{a,a+2} ,
\end{equation}
where
\begin{equation}
  \label{line_L_a_b}
  \begin{aligned}[b]
    E_{a,a+1}
    & =\left\{
      w_{a+3}=w_{a+4}=0
      \right\},
    \\
    E_{a,a+2}
    & =\left\{
      w_{a+2}+w_{a+3}=  w_{a+4}=0
      \right\}.
  \end{aligned}
\end{equation}
These  form a~configuration as in
Fig.~\ref{fig:pentagram}.

\subsubsection{$\mathbf{b}=(-2,-2,-2,-2,2)$}

We have
\begin{equation*}
  \Sing(\mathcal{Q})
  =\mathcal{S}_4 \cup \mathcal{S}_3 ,
\end{equation*}
and
five singular loci are given by
\begin{multline}
  \Sing(\widetilde{Q})
  =
  \left\{
    ( (-2,s,t,s,2), [1:-1:0:0:0] )
  \right\}
  \\
  \cup
    \left\{
    ( (-2,-2,s,t,s), [0:1:-1:0:0] )
  \right\}
  \cup
  \left\{
    ( (s,-2,-2,s,t), [0:0:1:-1:0] )
  \right\}
  \\
  \cup
    \left\{
      ( (t,s,-2,2,s), [0:0:0:1:-1] )
  \right\}
  \cup
  \left\{
    ( (s,t,-s,2,2), [1:0:0:0:1] )
  \right\}  .
\end{multline}
They intersect each other at distinct points, and
the jumping $\mathbb{P}^1$ fibers at the intersections
are as follows:
\begin{equation*}
  \begin{array}{cc}
    \toprule
    \mathbf{p}
    & \left. \Sing(\widetilde{\mathcal{Q}})
      \right|_{\mathbf{p}}
      \\
    \midrule
    (-2,-2,2,-2,2)
    &
      w_0+w_1+w_2=w_3=w_4=0
    \\
    (-2,-2,-2,-2,2)
    &
    w_0+w_1=w_2+w_3=w_4=0
    \\
    (-2,2,-2,2,2)
  &
    w_0+w_1=w_2=w_3+w_4=0
  \\
    (-2,2,2,2,2)
  &
    w_0+w_1-w_4=w_2=w_3=0
  \\
    (-2,-2,-2,-2,-2)
  &
    w_0=w_1+w_2+w_3=w_4=0
  \\
    (-2,-2,-2,2,-2)
  &
    w_0=w_1+w_2=w_3+w_4=0
  \\
    (-2,-2,2,2,2)
  &
    w_0-w_4=w_1+w_2=w_3=0
  \\
    (2,-2,-2,2,-2)
  &
    w_0=w_1=w_2+w_3+w_4=0
  \\
    (2,-2,-2,2,2)
  &
    w_0-w_4=w_1=w_2+w_3=0
  \\
    (2,2,-2,2,2)
  &
    w_0-w_3-w_4=w_1=w_2=0
    \\
    \bottomrule
  \end{array}
\end{equation*}
\subsection{Automorphism}

The mapping class  group of the  five-punctured sphere is generated by~$\sigma_a$
subject to the relations~\cite{Birman74,FarbMarg11Book}
\begin{equation}
  \label{Aut_Sigma05}
  \Aut (\Sigma_{0,5})=
  \left\langle
    \sigma_0, \sigma_1, \sigma_2, \sigma_3
    ~\middle|~
    \begin{matrix}
      \text{$
      \sigma_a  \, \sigma_{a+1} \, \sigma_a =
      \sigma_{a+1} \,  \sigma_{a} \, \sigma_{a+1}
      $ for $a=0,1,2$,}
      \\
      \text{$
      \sigma_a \, \sigma_b = \sigma_b \, \sigma_a
      $
      for $|a-b|>1$,}
      \\
      \sigma_0  \, \sigma_1 \, \sigma_2 \, \sigma_3
      \sigma_3 \, \sigma_2 \, \sigma_1 \, \sigma_0
      =
      \left(\sigma_0 \,  \sigma_1 \,  \sigma_2 \, \sigma_3\right)^5
      =1
    \end{matrix}
  \right\rangle .
\end{equation}
Each $\sigma_a$ is a half-twist exchanging  the $a$-th and the $(a+1)$-th
punctures
as in Fig.~\ref{fig:braid},
and it
swaps ${b}_a$ with
${b}_{a+1}$.

The actions are as follows
(see also~\cite{CalliMazzo18a}).
\begin{lemma}
The actions of $\sigma_a$ on the simple closed curves $\mathbf{c}$~\eqref{20-curve}
are given  by
\begin{subequations}
  \label{braid_on_curve}
  \begin{align}
    \sigma_0(\mathbf{c})
    &=
      \begin{bmatrix}
        {p}_{01}
        &{p}_{30}
        &      p_{02} \, p_{23}-p_{30}-b_0 \, b_3 - b_2\, p_{230}
        & b_1
        \\
        {p}_{02}
        &{p}_{24}
        &       p_{013}
        & b_0
        \\
        {p}_{23}
        &       {p}_{01} \, {p}_{30} - {p}_{13} - {b}_1\,  {b}_3
          - {b}_0 \, {p}_{013}
        & p_{402}
        & b_2
        \\
        {p}_{34}
        &      {p}_{40}
        &p_{230} \, p_{01} - p_{40} - b_0 \,  b_4 - b_1 \, p_{23}
        & b_3
        \\
        {p}_{230}
        &      {p}_{01} \, {p}_{02} - {p}_{12}   -{b}_1  \, {b}_2 - {b}_0 \,{p}_{34}
        &       p_{12}
        & b_4
      \end{bmatrix}
      ,
    \\
    \sigma_1(\mathbf{c})
    & =
      \begin{bmatrix}
        {p}_{341}
        & {p}_{12} \, {p}_{13} - {p}_{23}  -{b}_2 \, {b}_3 - {b}_1 \, {p}_{40}
        &     p_{23}
        & b_0
        \\
        {p}_{12}
        & {p}_{41}
        & p_{13} \, p_{34} - p_{41} - b_1 \, b_4 - b_3 \, p_{341}
        & b_2
        \\
        {p}_{13}
        & {p}_{30}
        &     p_{124}
        & b_1
        \\
        {p}_{34}
        &  {p}_{41} \, {p}_{12} - {p}_{24}   -{b}_2 \,  {b}_4 -{b}_1 \, {p}_{124}  
        &  p_{013}
        & b_3
        \\
        {p}_{40}
        &    {p}_{01}
        &     p_{341} \, p_{12} - p_{01} - b_0 \, b_1 - b_2 \, p_{34}
        & b_4
      \end{bmatrix}
      ,
    \\
    \sigma_2(\mathbf{c})
    & =
      \begin{bmatrix}
        {p}_{01}
        &     {p}_{12} 
        &     p_{402} \, p_{23} - p_{12} - b_1 \,  b_2 - b_3 \, p_{40}
        & {b}_0 
        \\
        {p}_{402}
        &  {p}_{24} \, {p}_{23} - {p}_{34}  -{b}_3 \, {b}_4 -{b}_2 \, {p}_{01} 
        &     p_{34}
        & {b}_1
        \\
        {p}_{23}
        &     {p}_{02}
        &     p_{124} \,  p_{23} - p_{341} - b_0 \,  b_2 - b_3 \,  p_{41} 
        & {b}_3
        \\
        {p}_{24}
        &     {p}_{41}
        &     p_{230}
        & {b}_2
        \\
        {p}_{40}
        & {p}_{02} \, {p}_{23} - {p}_{30}   -{b}_0  \, {b}_3 -{b}_2 \, {p}_{230} 
        &     p_{124}
        & {b}_4
      \end{bmatrix} ,
    \\
    \sigma_3(\mathbf{c})
    & =
      \begin{bmatrix}
        {p}_{01}
        & {p}_{13} \, {p}_{34} - {p}_{41}  -{b}_1 \, {b}_4 -{b}_3 \, {p}_{341} 
        & p_{230}
        & {b}_0 
        \\
        {p}_{12}
        & {p}_{23}
        & p_{013}  \, p_{34} - p_{23} - b_2 \,  b_3 - b_4 \, p_{01}
        & {b}_1 
        \\
        {p}_{013}
        & {p}_{30} \, {p}_{34} - {p}_{40}  -{b}_0 \, {b}_4 -{b}_3 \, {p}_{12} 
        & p_{40}
        & {b}_2
        \\
        {p}_{34}
        & {p}_{13}
        & p_{230} \, p_{34} - p_{402} - b_1 \,  b_3 - b_4 \,  p_{02}
        &{b}_4  
        \\
        {p}_{30}
        & {p}_{02}
        & p_{341}
        & {b}_3
  \end{bmatrix}
    .
  \end{align}
\end{subequations}
\end{lemma}
\section{Cluster Algebra}
\label{sec:cluster}
\subsection{Cluster Algebra and  Ideal Triangulation}
We briefly recall the definition of a cluster algebra~\cite{FominZelev02a}.
Let $Q$ be a finite quiver without  loops.
Each vertex of $Q$ is assigned a cluster variable $y_i$
for $i=1,\dots,|Q|$,
and we define a skew-symmetric matrix $\mathbf{B}=(b_{ij})$ by
\begin{equation}
  \label{define_b_ij}
  b_{ij}
  = \#(\text{arrows from $\overset{i}{\bullet}$ to  $\overset{j}{\bullet}$})
  - \#(\text{arrows from $\overset{j}{\bullet}$ to  $\overset{i}{\bullet}$})
\end{equation}
A seed is a pair $(\mathbf{y},\mathbf{B})$,
and 
a mutation $\mu_k$ at the vertex $k$ is an involutive  map,
\begin{equation}
  \label{define_mu}
  \mu_k(\mathbf{y},\mathbf{B})
  =
  (\widetilde{\mathbf{y}}, \widetilde{\mathbf{B}}) .
\end{equation}
Here
\begin{equation}
  \label{mutation_y_i}
  \widetilde{y}_i
  =
  \begin{cases}
    y_k^{-1}, & \text{for $i=k$,}
    \\
    y_i \left( 1+ y_k^{-1} \right)^{-b_{ki}},
              & \text{for $i\neq k$
                and $b_{ki}\geq 0$,}
    \\
    y_i \left( 1+ y_k \right)^{-b_{ki}},
              & \text{for $i\neq k$
                and $b_{ki}\leq 0$,}
  \end{cases}
\end{equation}
and an exchange matrix $\widetilde{\mathbf{B}}=(\widetilde{b}_{ij})$ is given by
\begin{equation}
  \label{mutation_b_ij}
  \widetilde{b}_{ij}
  =
  \begin{cases}
    -b_{ij}, & \text{for $i=k$ or $j=k$,}
    \\
    b_{ij} + \frac{|b_{ik}|b_{kj}+b_{ik}|b_{kj}|}{2},
             & \text{otherwise.}
  \end{cases}
\end{equation}

The cluster algebra has a Poisson structure defined by the exchange
matrix $\mathbf{B}$
(see, \emph{e.g.},~\cite{GekhShapVain10Book}).
The cluster variables $y_i$ serve
as log-canonical coordinates of the
Poisson structure, and the Poisson bracket is defined by
\begin{equation}
  \label{Poisson_y-y-B}
  \left\{ y_i ~,~ y_j \right\}
  = b_{ij}  \, y_i \,  y_j.
\end{equation}
Cluster mutations~\eqref{define_mu} preserve this Poisson
structure.

\begin{figure}[htbp]
  \centering
  \begin{minipage}{0.4\linewidth}
    \centering
    \includegraphics[scale=.9]{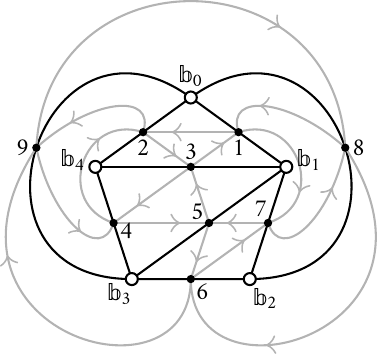}
  \end{minipage}
  \qquad
  \begin{minipage}{0.4\linewidth}
    \centering
    \includegraphics[scale=.9]{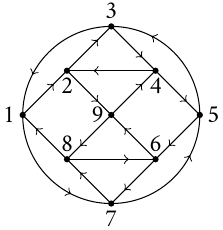}
  \end{minipage}
  \caption{An ideal triangulation of $\Sigma_{0,5}$ is depicted in the
    left panel.
    Each puncture is denoted as a white circle,
    and it is labeled by the boundary curve $\mathbb{b}_a$.
    The vertices of
    the ideal triangulations are punctures, and
    we assign a vertex~$\bullet$
    of a quiver in each edge of the ideal triangles.
    The resulting quiver $Q_{\Sigma_{0,5}}$ is rewritten in the right panel.}
  \label{fig:ideal_triangulation}
\end{figure}

We consider the cluster algebra
associated with
the ideal triangulation of~$\Sigma_{0,5}$
shown in
Fig.~\ref{fig:ideal_triangulation}.
By assigning  counterclockwise arrows for each triangle,
the ideal triangulation gives the quiver $Q_{\Sigma_{0,5}}$ in Fig.~\ref{fig:ideal_triangulation},
whose exchange matrix is
\begin{equation}
  \label{eq:7}
  \mathbf{B}
  =
  \begin{pmatrix}
    &1 &-1&&&&1&-1&\\
    -1&&1&-1&&&&&1\\
    1&-1&&1&-1&&&&\\
    &1&-1&&1&&&&-1\\
    &&1&-1&&1&-1&&\\
    &&&&-1&&1&-1&1\\
    -1&&&&1&-1&&1&\\
    1&&&&&1&-1&&-1\\
    &-1&&1&&-1&&1&
  \end{pmatrix}.
\end{equation}

It is known that the simple closed curves can be written as Laurent
polynomials in the square roots of
the cluster variables.
Following the method in~\cite{BonaWong16b,TQLe16a}
(see also~\cite{KHikami17a} for $\Sigma_{0,4}$ and $\Sigma_{1,1}$,
and~\cite{KHikami24b} for $\Sigma_{1,2}$),
the boundary circles are given by
\begin{equation}
  \label{boundary_cluster}
  \begin{aligned}[b]
    {b}_0
    & \mapsto
      \sqrt{y_1y_2 y_8 y_9}
      +\frac{1}{  \sqrt{y_1y_2 y_8 y_9}} ,
    \\
    {b}_1
    &\mapsto
      \sqrt{y_1y_3 y_5 y_7}
      +
      \frac{1}{\sqrt{y_1y_3 y_5 y_7}} ,
    \\
    {b}_2
    & \mapsto
      \sqrt{y_6y_7 y_8}+
      \frac{1}{\sqrt{y_6y_7 y_8}} ,
    \\
    {b}_3
    &\mapsto
      \sqrt{y_4y_5 y_6 y_9}
      +\frac{1}{\sqrt{y_4y_5 y_6 y_9}} ,
    \\
    {b}_4
    &\mapsto
      \sqrt{y_2 y_3 y_4}
      +\frac{1}{\sqrt{y_2 y_3 y_4}} .
  \end{aligned}
\end{equation}
In the same manner,
we find that
the simple closed curves~${p}_{a,a+1}$ are given by
\begin{equation}
  \label{p_a_a+1_cluster}
  \begin{aligned}[b]
    {p}_{01}
    &\mapsto   \alpha_6(y_8,y_3,y_5,y_2,y_7,y_9) ,
    \\
    {p}_{12}
    & \mapsto \alpha_5(y_1,y_6,y_3,y_5,y_8) ,
    \\
    {p}_{23}
    & \mapsto \alpha_5(y_5,y_8,y_4,y_9,y_7),
    \\
    {p}_{34}
    & \mapsto \alpha_5(y_9,y_3,y_6,y_5,y_2),
    \\
    {p}_{40}
    & \mapsto \alpha_5(y_1,y_4,y_8,y_9,y_3),
  \end{aligned}
\end{equation}
where
\begin{align}
  &
    \begin{multlined}[b][.9\textwidth]
        \alpha_5(y_1,y_2,y_3,y_4,y_5)
  =
    \sqrt{y_1 y_2 y_3 y_4 y_5}
    \left(1+\frac{1}{y_1} \right) \left(1+\frac{1}{y_2} \right)
    \\
    +\sqrt{\frac{y_3 y_4}{y_1 y_2 y_5}} +
    \sqrt{\frac{y_2 y_4 y_5}{y_1 y_3}}
    +
    \frac{1}{\sqrt{y_1 y_2 y_3 y_4 y_5}}
    (1+y_4)(1+y_5) ,
  \end{multlined}
  \\
  &
    \begin{multlined}[b][.9\textwidth]
      \alpha_6(y_1,y_2,y_3,y_4,y_5,y_6)
      =
      \sqrt{y_1 y_2 y_3 y_4 y_5 y_6}
      \left(1+\frac{1}{y_1}+\frac{1}{y_1 y_6} \right)
      \left(1+\frac{1}{y_2} \right)
      \\
      +\sqrt{\frac{y_3 y_5}{y_1 y_2 y_4 y_6}} +
      \sqrt{\frac{y_1y_4 y_5 y_6}{y_2 y_3}}
      +
      \frac{1}{\sqrt{y_1 y_2 y_3 y_4 y_5 y_6}}
      (1+y_4+y_4 y_6)(1+y_5) .
  \end{multlined}
\end{align}
The simple closed curves ${p}_{a,a+2}$ and
${p}_{a,a+1,a+3}$
are also given by
\begin{equation}
  \label{p_a_a+2_cluster}
  \begin{aligned}[b]
    {p}_{341}
    & \mapsto \alpha_5(y_9,y_7,y_2,y_1,y_6) ,
    \\
    {p}_{13}
    & \mapsto \alpha_6(y_7,y_4,y_9,y_3,y_6, y_1) ,
    \\
    {p}_{24}
    & \mapsto
      \alpha_7(y_2,y_3,y_4,y_5,y_6,y_7, y_8) ,
    \\
    {p}_{124}
    & \mapsto \alpha_6(y_6,y_2,y_1,y_4,y_8, y_5) ,
    \\
    {p}_{41}
    & \mapsto \alpha_5(y_5,y_2,y_7,y_1,y_4) ,
  \end{aligned}
\end{equation}
where
\begin{multline}
  \alpha_7(y_1,y_2,y_3,y_4,y_5,y_6,y_7)
  =
  \sqrt{\frac{y_3}{y_1 y_2  y_5 y_6 y_7}}
  (1+y_1) \left(1+\frac{1}{y_4} \right)
  (1+y_4 y_5 y_6 y_7)
  \\
  + \sqrt{\frac{y_6}{y_1 y_2  y_3 y_5 y_7}}
  \left(1+\frac{1}{y_4} \right)(1+y_7)
  (1+y_1 y_2 y_3 y_4)
  +
  \sqrt{\frac{y_3 y_6}{y_1 y_2 y_5 y_7}}
  (1+y_1)   \frac{(1+y_4)^2}{y_4}  (1+y_7)
  \\
  + y_4 \sqrt{y_1 y_2 y_3 y_5 y_6 y_7}
  + \frac{1}{ y_4 \sqrt{y_1 y_2 y_3 y_5 y_6 y_7}} .
\end{multline}
We omit some of the explicit expressions for ${p}_{a,a+2}$ and
${p}_{a,a+1,a+3}$, as they  can be obtained from~\eqref{pp-1}
using~\eqref{boundary_cluster},~\eqref{p_a_a+1_cluster} and~\eqref{p_a_a+2_cluster}.

We can check by \texttt{Mathematica}  that these maps
fulfill~\eqref{W_for_2-p},~\eqref{our-M-eq},~\eqref{rank-Q},~\eqref{vector_kernel},
and~\eqref{adj-Q_v_v}:
\begin{prop}
  The simple closed curves $\mathbf{c}$~\eqref{20-curve}
  on~$\Sigma_{0,5}$
  are written
  in terms of the cluster $y$-variables
  as~\eqref{boundary_cluster},~\eqref{p_a_a+1_cluster}, and~\eqref{p_a_a+2_cluster}.
\end{prop}

\subsection{Goldman Poisson Algebra}

The skein algebra $\Sk_A(\Sigma)$ gives the Poisson structure for the
simple closed curves on the surface
in the limit $A \to \pm 1$.
In the case of
$\Sigma_{0,4}$,
the curves in Fig.~\ref{fig:curves_Sigma04}
satisfy the  Poisson
algebra
(see, \emph{e.g.}, \cite{NekrRoslShat11a,ChekhMazzo10a,EtingGinzb10a}).
We recall that
\begin{equation}
  \label{eq:3}
  \begin{aligned}[b]
    \left\{  {p}_{01} ~,~ {p}_{12} \right\}
    & =
      -{p}_{12} \, {p}_{01} + 2 \, {p}_{20}
      +
      {b}_{1} \, {b}_{3} + {b}_{0} \, {b}_{2}
      =
      \frac{\partial W}{\partial {p}_{20}}
      ,
    \\
    \left\{  {p}_{20} ~,~ {p}_{01} \right\}
    & =
      - {p}_{01} \, {p}_{20} + 2\, {p}_{12}
      +
      {b}_{1} \, {b}_{2} + {b}_{0} \, {b}_{3}
      =
        \frac{\partial W}{\partial {p}_{12}}
      ,
    \\
    \left\{  {p}_{12} ~,~ {p}_{20} \right\}
    & =
      -{p}_{20} \, {p}_{12} + 2 \, {p}_{01}
      +
      {b}_{0}  \, {b}_{1} + {b}_{2} \, {b}_{3}
      =
        \frac{\partial W}{\partial {p}_{01}}
      ,
  \end{aligned}
\end{equation}
where $W=  W_{{b}_0,{b}_1,{b}_2,{b}_3}
  ({p}_{01}, {p}_{12}, {p}_{20})
$ as in~\eqref{Fricke-Klein}.
We denote this set of Poisson relations as
$\Pois_{{b}_0,{b}_1,{b}_2,{b}_3}
({p}_{01}, {p}_{12}, {p}_{20}) $.
The boundary curves ${b}_a$ are  central in the Poisson algebra,
and
their Poisson brackets with all other curves vanish.

For~$\Sigma_{0,5}$,
the boundary curves~$b_a$ are central as well.
The  Poisson algebra for the curves in
Fig.~\ref{fig:curves_Sigma05} is induced from~\eqref{Poisson_y-y-B}
as a classical limit~$A\to 1$ of the commutation relations in
$\Sk_A(\Sigma_{0,5})$~\cite{CookeLacab22a}.

\begin{prop}
  The Goldman Poisson relations of the simple closed curves
  $\mathbf{c}$~\eqref{20-curve}
  are given by
\begin{equation}
  \label{eq:13}
  \begin{gathered}[b]
  \Pois_{{b}_a,{b}_{a+1},{b}_{a+2},{p}_{a+3,a+4}}
  ({p}_{a,a+1}, {p}_{a+1,a+2}, {p}_{a,a+2}),
  \\
  \Pois_{{b}_{a+1},{b}_{a},{b}_{a+2},{p}_{a+3,a+4}}
  ({p}_{a,a+1}, {p}_{a+3,a+4,a+1}, {p}_{a,a+1}),
  \\
  \Pois_{{b}_{a},{b}_{a+2},{b}_{a+4},{p}_{a+4,a,a+2}}
  ({p}_{a,a+2}, {p}_{a+2,a+4}, {p}_{a+4,a}),
  \\
  \Pois_{{b}_{a+2},{b}_{a},{b}_{a+4},{p}_{a+1,a+3}}
  ({p}_{a+3,a+4,a+1}, {p}_{a+4,a},
  {p}_{a,a+1,a+3}),
  \\
  \begin{multlined}[b][.88\textwidth]
    \left\{ {p}_{a,a+2} ~,~ {p}_{a+1,a+3} \right\}
    =
    \frac{1}{2}
    \left(
      2 {p}_{a+1,a+2}+ {b}_{a+1} {b}_{a+2}
    \right)
    \left(
      2 {p}_{a+3,a}+ {b}_{a+3} {b}_{a}
    \right)
    \\
    -
    \frac{1}{2}
    \left(
      2 {p}_{a,a+1}+ {b}_{a} {b}_{a+1}
    \right)
    \left(
      2 {p}_{a+2,a+3}+ {b}_{a+2} {b}_{a+3}
    \right) ,
  \end{multlined}
  \\
    \begin{multlined}[b][.88\textwidth]
      \left\{ {p}_{a+3,a+4,a+1} ~,~ {p}_{a+4,a,a+2} \right\}
    =
    \frac{1}{2}
    \left(
      2 {p}_{a,a+1}+ {b}_{a} {b}_{a+1}
    \right)
    \left(
      2 {p}_{a+2,a+3}+ {b}_{a+2} {b}_{a+3}
    \right)
    \\
    -
    \frac{1}{2}
    \left(
      2 {p}_{a+1,a+2}+ {b}_{a+1} {b}_{a+2}
    \right)
    \left(
      2 {p}_{a+1,a+2,a+4}+ {b}_{a} {b}_{a+3}
    \right) .
  \end{multlined}
\end{gathered}
\end{equation}
Two simple closed curves without intersection have vanishing Poisson brackets.
\end{prop}

\subsection{Automorphisms}
\begin{figure}[htbp]
  \centering
  \begin{tikzpicture}[
    img/.style={inner sep=0pt, outer sep=0pt},
      mapsto/.style={
        |->,
        >={To[length=4pt, width=4pt]}, %
        line width=0.6pt, 
        shorten <=5pt, 
        shorten >=5pt
    }
    ]

    \def\R{55mm}

    \foreach \file [count=\i] in {pentagon_triangle-Init, pentagon_triangle-sigma0, pentagon_triangle-sigma1, pentagon_triangle-sigma2, pentagon_triangle-sigma3} {

        \pgfmathsetmacro{\angle}{90 - (\i - 1) * 72}

        \node[img] (p\i) at (\angle:\R) {
            \includegraphics[scale=.8]{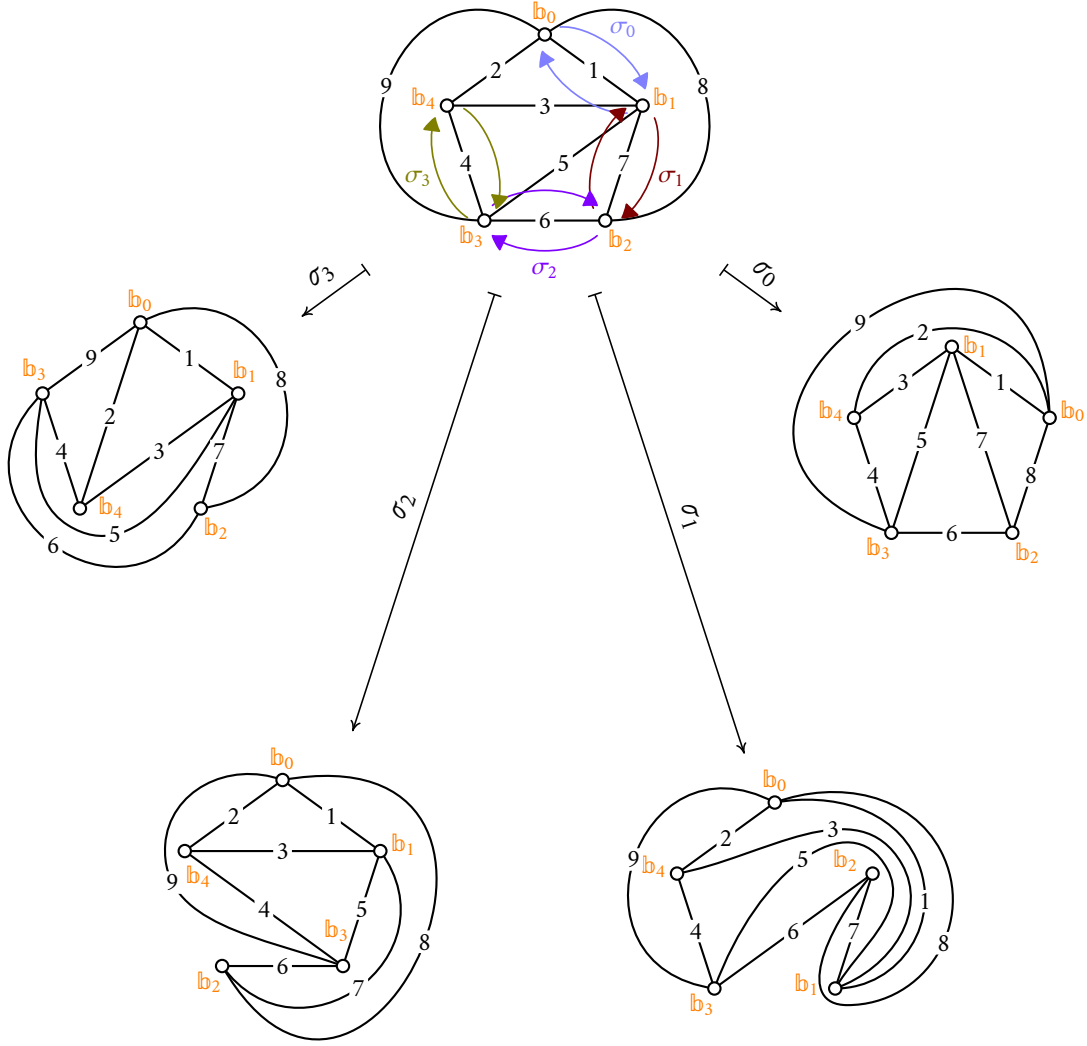}
        };

    }

    \draw[mapsto] (p1) -- node[midway, auto, sloped] {$\sigma_0$}  (p2);
    \draw[mapsto] (p1) -- node[midway, auto, sloped] {$\sigma_1$}  (p3);
    \draw[mapsto] (p1) -- node[midway, auto, sloped] {$\sigma_2$}  (p4);
    \draw[mapsto] (p1) -- node[midway, auto, sloped] {$\sigma_3$}  (p5);
\end{tikzpicture}
  \caption{The braid operator $\sigma_a$ acts as an exchange of the $a$-th
    and the $(a+1)$-th punctures.
    The labels of edge denote those of the
    cluster $y$-variables in Fig.~\ref{fig:ideal_triangulation}.}
\label{fig:braid}
\end{figure}

\begin{figure}[htbp]
  \centering
    \begin{tikzpicture}[
    img/.style={inner sep=0pt, outer sep=0pt},
    mapsto/.style={{Bar[width=4pt, line width=0.6pt]}-{To[length=4pt, width=4pt]}, line width=0.6pt, shorten <=5pt, shorten >=5pt}
    ]
    \node[img] (img1) {\includegraphics[scale=.8]{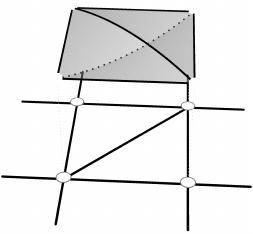}};
  \node[img, right=15mm of img1] (img2) {\includegraphics[scale=.8]{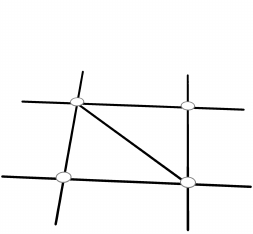}};
  \draw[mapsto] (img1.east) -- (img2.west) node[midway,below]{};
\end{tikzpicture}
\caption{Depicted is the flip of an ideal triangulation.
  We can interpret the flip 
  as an attachment of an ideal hyperbolic tetrahedron to the surface.
  See~\cite{KHikamiRInoue14a}.
}
  \label{fig:flip}
\end{figure}

With the simple closed curves  expressed in terms of cluster $y$-variables,
the mapping class group~\eqref{Aut_Sigma05} admits a realization  in terms of
cluster mutations
(see~\cite{KHikami17a,KHikami24b} for~$\Sigma_{0,4}$ and $\Sigma_{1,2}$).
We recall the ideal triangulations in
Fig.~\ref{fig:ideal_triangulation}.
Observe  that $\sigma_a$ acts as a braiding operator as in
Fig.~\ref{fig:braid},
and that we obtain  a different ideal triangulation from the initial one.
Since cluster   mutations correspond to  flips of ideal
triangulations~\cite{FominZelev02a} as in Fig.~\ref{fig:flip},
each braid generator is  realized as a composition of cluster mutations:
\begin{equation}
  \label{braid_mutation}
  \begin{aligned}[b]
    \sigma_0
    & = s_{29} \,  s_{28} \,  s_{27} \,  s_{25} \, s_{23} \, \mu_{9,8,5,3} ,
    \\
    \sigma_1
    & = s_{18} \,  s_{16} \,  s_{15} \,  s_{13} \, \mu_{3,1,6} ,
    \\
    \sigma_2
    & = s_{45}  \, s_{47} \,  s_{48} \, s_{49}  \, \mu_{8,4,5} ,
    \\
    \sigma_3
    & = s_{23} \,  s_{25} \, s_{26}\,  s_{29}  \,  \mu_{6,9,3} .
  \end{aligned}
\end{equation}
Here $s_{i,j}$ denotes the transposition of the labels of~$i$ and~$j$,
and
$\mu_{i,j,\dots,k}$ means $\mu_i \mu_j \cdots \mu_k$.

\begin{lemma}
  The generators of the automorphism~\eqref{Aut_Sigma05}
  act
  on the cluster $y$-variables
  $\mathbf{y}=(y_1,\dots,y_9)^\top$ as follows;
  \begin{subequations}
    \begin{gather}
      \sigma_0(\mathbf{y}) =
      \begin{pmatrix}
        \frac{y_1 y_3 y_5 y_8 y_9}{(1 + y_5 + y_3 y_5) (1 + y_9 + y_8 y_9)}
        \\
        \frac{1}{y_9 (1+  y_8) }
        \\
        y_2 (1 + y_3) (1 + y_9 + y_8 y_9)\\
        \frac{y_3 y_4 (1 + y_8) y_9}{(1 +    y_3) (1 + y_9 + y_8 y_9)}
        \\
        \frac{1 + y_5 +  y_3 y_5}{y_3}
        \\
        \frac{(1 + y_3) y_5 y_6 y_8}{(1 + y_5 + y_3 y_5) (1 + y_8)} \\
        \frac{1}{y_5(1 +       y_3 )}
        \\
        (1 + y_5 + y_3 y_5) y_7 (1 + y_8)
        \\
        \frac{1 + y_9 + y_8 y_9}{y_8}
      \end{pmatrix},
      \quad
      \sigma_1(\mathbf{y}) =
      \begin{pmatrix}
        (1 + y_1) (1 + y_6) y_8      \\
        \frac{y_1 y_2}{1 + y_1} \\
        \frac{ 1 + y_3 + y_1 y_3}{y_1} \\
        \frac{(1 + y_1) y_3 y_4}{1 + y_3 + y_1 y_3}\\
        \frac{1}{ y_3(1 + y_1)}\\
        (1 + y_3 + y_1 y_3) y_5 (1 + y_6)\\
        \frac{ y_1 y_3 y_6 y_7}{(1 + y_3 + y_1 y_3) (1 + y_6)}\\
        \frac{1}{y_6}\\
        \frac{y_6 y_9}{1 + y_6}      
      \end{pmatrix}
      ,
    \\
    \sigma_2(\mathbf{y}) =
    \begin{pmatrix}
      \frac{y_1 y_8}{1 + y_8} \\
      \frac{y_2 y_4 (1 + y_5)}{1 + y_4 +  y_4 y_5 } \\
      \frac{y_3 y_5}{1 + y_5} \\
      \frac{1 + y_4 +  y_4 y_5}{y_5} \\
      (1 + y_5) y_7 (1 + y_8) \\
      \frac{y_4 y_5 y_6 y_8}{(1 + y_4 + y_4 y_5) (1 +      y_8)} \\
      \frac{1}{y_8} \\
      (1 + y_4 + y_4 y_5) (1 + y_8) y_9 \\
      \frac{1}{y_4(1 +  y_5)      }
    \end{pmatrix},
    \quad
    \sigma_3(\mathbf{y}) =
    \begin{pmatrix}
      \frac{y_1 y_3}{1 + y_3} \\
      \frac{1}{y_3} \\
      (1 + y_3) y_5 (1 + y_6 +    y_6 y_9) \\
      \frac{y_3 y_4 y_6 y_9}{(1 + y_3) (1 + y_6 + y_6 y_9)} \\
      \frac{1}{y_6(1 +   y_9)} \\
      \frac{1 + y_6 + y_6 y_9}{y_9} \\
      \frac{y_6 y_7 (1 + y_9)}{1 + y_6 +  y_6 y_9} \\
      \frac{y_8 y_9}{1 + y_9} \\
      y_2 (1 + y_3) (1 + y_9)      
    \end{pmatrix} .
  \end{gather}
\end{subequations}
\end{lemma}

It is straightforward to verify by \texttt{Mathematica}
the following statements.

\begin{prop}
  The mutations~\eqref{braid_mutation}
  generate  $\Aut(\Sigma_{0,5})$~\eqref{Aut_Sigma05}.
\end{prop}

\begin{lemma}
  The mutations~\eqref{braid_mutation}
  induce~\eqref{braid_on_curve}.
\end{lemma}

\section{Concluding Remarks}
We have derived the determinantal quintic  hypersurface~\eqref{det-Q} associated with the
character variety of the five-punctured sphere using the cluster
$y$-variables for the simple closed curves on the surface.
The explicit form~\eqref{det-Q-explicit} is a generalization of
Simpson's hypersurface~\cite{Simpson18a} to arbitrary boundary
conjugacy classes.
The method developed  here may be extended to the character variety of
$\Sigma_{g,n}$ for general $(g,n)$,
and 
to the
construction of refined  quantum invariants via
the cluster algebra and the double affine Hecke
algebra,
following the programme of~\cite{KHikami25a,KHikami24b,KHikami24a,KHikami19a}.

\section*{Acknowledgements}
The work of KH is supported in part by
JSPS KAKENHI
Grant Number JP23K22388, JP22H01117.

\section*{Use of AI}
ChatGPT and Claude were used for discussions.
The author  has also used them to rewrite the English and to produce the
\texttt{SageMath} and asymptote codes.
He takes full responsibility for all contents.

\appendix
\section{\texttt{SageMath}}

For a  specific value of $\mathbf{b}$,
say $\mathbf{b}=(0,0,0,0,0)$,
we may use the following code.
\begin{lstlisting}[style=sageinput]
R.<x,y,z,u,v> = QQ[]
Vars = [x,y,z,u,v]
M = matrix(R, [
  [2,0,-x,-u,0],
  [0,2,0,-y,-v],
  [-x,0,2,0,-z],
  [-u,-y,0,2,0],
  [0,-v,-z,0,2]
])
f = M.det()
derivatives = [f.derivative(q) for q in R.gens()]
I = R.ideal([f] + derivatives)
I.dimension()
@\sageoutstyle{2}@
I4 = R.ideal(M.minors(4))
I == I4
@\sageoutstyle{True}@
if not I==I4:
  I.quotient(I4)

I3 = R.ideal(M.minors(3))
J = I + I3
J.dimension()
@\sageoutstyle{-1}@
if not J.dimension()==-1:
  for pt in I3.variety(QQ):
    print(" ", pt)
\end{lstlisting}

%
%

\begin{thebibliography}{10}
\providecommand{\url}[1]{\texttt{#1}}
\providecommand{\urlprefix}{URL }
\providecommand{\eprint}[2][]{\url{#2}}

\bibitem{Beauv00a}
A.~Beauville, \emph{Determinantal hypersurfaces},
  \href{http://dx.doi.org/10.1307/mmj/1030132707}{Michigan Math. J.}
  \textbf{48}, 39--64 (2000).

\bibitem{Birman74}
J.~S. Birman, \emph{Braids, Links, and Mapping Class Groups}, vol.~83 of
  \emph{Annals of Math. Stud.}, Princeton Univ. Press, Princeton, 1974.

\bibitem{BonaWong16b}
F.~Bonahon and H.~Wong, \emph{Representations of the {Kauffman} bracket skein
  algebra {I}: invariants and miraculous cancellations},
  \href{http://dx.doi.org/10.1007/s00222-015-0611-y}{Invent. Math.}
  \textbf{204}, 195--243 (2016).

\bibitem{Bullo97a}
D.~Bullock, \emph{Rings of {$SL_2(\mathbb{C})$}-characters and the {Kauffman}
  bracket skein module},
  \href{http://dx.doi.org/10.1007/s000140050032}{Comment. Math. Helv.}
  \textbf{72}, 521--542 (1997).

\bibitem{CalliMazzo18a}
P.~Calligaris and M.~Mazzocco, \emph{Finite orbits of the pure braid group on
  the monodromy of the 2-variable {Garnier} system},
  \href{http://dx.doi.org/10.1093/integr/xyy005}{J. Integrable Systems}
  \textbf{3}, 1--35 (2018).

\bibitem{ChekhMazzo10a}
L.~Chekhov and M.~Mazzocco, \emph{Shear coordinate description of the quantized
  versal unfolding of a {$D_4$} singularity},
  \href{http://dx.doi.org/10.1088/1751-8113/43/44/442002}{J. Phys. A: Math.
  Theor.} \textbf{43}, 442002 (2010), 13 pages.

\bibitem{CookeLacab22a}
J.~Cooke and A.~Lacabanne, \emph{Higher rank {Askey--Wilson} algebras as skein
  algebras}, \href{http://dx.doi.org/10.5802/aif.3729}{Ann. Inst. Fourier}
  \textbf{76}, 1055--1117 (2025).

\bibitem{Dolga12Book}
I.~V. Dolgachev, \emph{Classical Algebraic Geometry: a modern view}, Cambridge
  Univ. Press, Cambridge, 2012.

\bibitem{EtingGinzb10a}
P.~Etingof and V.~Ginzburg, \emph{Noncommutative del {Pezzo} surfaces and
  {Calabi}--{Yau} algebras}, \href{http://dx.doi.org/10.4171/JEMS/235}{J. Euro.
  Math. Soc.} \textbf{12}, 1371--1416 (2010).

\bibitem{TehraFrohm26a}
M.~Farajzadeh-Tehrani and C.~Frohman, \emph{On compactifications of the
  {$SL(2,\mathbb{C}$)} character varieties of punctured surfaces},
  \href{http://dx.doi.org/10.1090/tran/9564}{Trans. Amer. Math. Soc.}
  \textbf{379}, 627--679 (2026).

\bibitem{FarbMarg11Book}
B.~Farb and D.~Margalit, \emph{A Primer on Mapping Class Groups}, vol.~49 of
  \emph{Princeton Math. Series}, Princeton Univ. Press, Princeton, 2011.

\bibitem{FassaLoray25a}
T.~Fassarella and F.~Loray, \emph{Moduli of {Higgs} bundles over the five
  punctured sphere}, \href{http://dx.doi.org/10.4310/AJM.260122041319}{Asian J.
  Math.} \textbf{29}, 733--759 (2025).

\bibitem{FominZelev02a}
S.~Fomin and A.~Zelevinsky, \emph{Cluster algebras {I}. foundations},
  \href{http://dx.doi.org/10.1090/S0894-0347-01-00385-X}{J. Amer. Math. Soc.}
  \textbf{15}, 497--529 (2002).

\bibitem{FrickeKleinBook1897}
R.~Fricke and F.~Klein, \emph{Vorlesungen {\"uber} die Theorie der Automorphen
  Functionen, Band 1: Die gruppentheoretischen Grundlagen}, Druck und Verlag
  von B. G. Teubner, Leipzig, 1897.

\bibitem{GekhShapVain10Book}
M.~Gekhtman, M.~Shapiro, and A.~Vainshtein, \emph{Cluster Algebra and Poisson
  Geometry}, vol. 167 of \emph{Mathematical Surveys and Monographs}, Amer.
  Math. Soc., Providence, 2010.

\bibitem{Girand16a}
A.~Girand, \emph{A new two-parameter family of isomonodromic deformations over
  the five punctured sphere}, \href{http://dx.doi.org/10.24033/bsmf.2716}{Bull.
  Soc. Math. France} \textbf{144}, 339--368 (2016).

\bibitem{Goldman09a}
W.~M. Goldman, \emph{Trace coordinates on {Fricke} spaces of some simple
  hyperbolic surfaces}, in A.~Papadopoulos, ed., \emph{Handbook of
  Teichm{\"u}ller Theory}, vol.~2, chap.~15, pp. 611--684, Euro. Math. Soc.,
  2009.

\bibitem{KHikami19a}
K.~Hikami, \emph{{DAHA} and skein algebra of surfaces: double-torus knots},
  \href{http://dx.doi.org/10.1007/s11005-019-01189-5}{Lett. Math. Phys.}
  \textbf{109}, 2305--2358 (2019).

\bibitem{KHikami17a}
---{}---{}---, \emph{Note on character varieties and cluster algebras},
  \href{http://dx.doi.org/10.3842/SIGMA.2019.003}{SIGMA} \textbf{15}, 003
  (2019), 32 pages.

\bibitem{KHikami24a}
---{}---{}---, \emph{Generalized double affine {Hecke} algebra for double
  torus}, \href{http://dx.doi.org/10.1007/s11005-024-01848-2}{Lett. Math.
  Phys.} \textbf{114}, 102 (2024), 16 pages.

\bibitem{KHikami25a}
---{}---{}---, \emph{{Askey--Wilson} polynomial and genus-two skein module},
  \href{http://dx.doi.org/10.1088/1751-8121/ae144a}{J. Phys. A: Math. Gen.}
  \textbf{58}, 445203 (2025), 12 pages.

\bibitem{KHikami24b}
---{}---{}---, \emph{A note on double affine {Hecke} algebra for skein algebra
  on twice-punctured torus},
  \href{http://dx.doi.org/10.1016/j.geomphys.2024.105408}{J. Geom. Phys.}
  \textbf{209}, 105408 (2025), 13 pages.

\bibitem{KHikamiRInoue14a}
K.~Hikami and R.~Inoue, \emph{Braiding operator via quantum cluster algebra},
  \href{http://dx.doi.org/10.1088/1751-8113/47/47/474006}{J. Phys. A: Math.
  Theor.} \textbf{47}, 47006 (2014), 21 pages.

\bibitem{RHorow72a}
R.~D. Horowitz, \emph{Characters of free groups represented in the two
  dimensional special linear group},
  \href{http://dx.doi.org/10.1002/cpa.3160250602}{Comm. Pure Appl. Math.}
  \textbf{25}, 635--649 (1972).

\bibitem{Iwasaki03a}
K.~Iwasaki, \emph{An area-preserving action of the modular group on cubic
  surfaces and the {Painlev{\'e}} {VI} equation},
  \href{http://dx.doi.org/10.1007/s00220-003-0940-3}{Commun. Math. Phys.}
  \textbf{242}, 185--219 (2003).

\bibitem{Komyo15a}
A.~Komyo, \emph{On compactifications of character varieties of $n$-punctured
  projective line}, \href{http://dx.doi.org/10.5802/aif.2965}{Ann. Inst.
  Fourier} \textbf{65}, 1493--1523 (2015).

\bibitem{TQLe16a}
T.~T.~Q. Le, \emph{Triangular decomposition of skein algebras},
  \href{http://dx.doi.org/10.4171/QT/115}{Quantum Topol.} \textbf{9}, 591--632
  (2018).

\bibitem{NekrRoslShat11a}
N.~A. Nekrasov, A.~A. Rosly, and S.~L. Shatashvili, \emph{{Darboux}
  coordinates, {Yang--Yang} functional, and gauge theory},
  \href{http://dx.doi.org/10.1016/j.nuclphysbps.2011.04.150}{Nucl. Phys. B:
  Proc. Suppl.} \textbf{216}, 69--93 (2011).

\bibitem{Piont06a}
J.~Piontkowski, \emph{Linear symmetric determinantal hypersurfaces},
  \href{http://dx.doi.org/10.1307/mmj/1144437441}{Michigan Math. J.}
  \textbf{54}, 117--155 (2006).

\bibitem{PrzytSikor19a}
J.~H. Przytycki and A.~S. Sikora, \emph{Skein algebra of surfaces},
  \href{http://dx.doi.org/10.1090/tran/7298}{Trans. Amer. Math. Soc.}
  \textbf{371}, 1309--1332 (2019).

\bibitem{Simpson16a}
C.~T. Simpson, \emph{The dual boundary complex of the {$SL_2$} character
  variety of a punctured sphere},
  \href{http://dx.doi.org/10.5802/afst.1496}{Ann. Fac. Sci. Toulouse Math.}
  \textbf{25}, 317--361 (2016).

\bibitem{Simpson18a}
---{}---{}---, \emph{An explicit view of the {Hitchin} fibration on the {Betti}
  side for {$\mathbb{P}^1$} minus five points}, in J.~E. Andersen, A.~Dancer,
  and O.~Garc{\'{i}}a-Prada, eds., \emph{Geometry and Physics: A Festschrift in
  Honour of {Nigel} {Hitchin}}, pp. 705--724, Oxford Univ. Press, Oxford, 2018.

\bibitem{Szabo23a}
S.~Szab{\'o}, \emph{Hitchin {WKB}-problem and {$P=W$} conjecture in lowest
  degree for rank 2 over the 5-punctured sphere},
  \href{http://dx.doi.org/10.1093/qmath/haac037}{Quarterly J. Math.}
  \textbf{74}, 687--746 (2023).

\bibitem{Vogt89a}
H.~Vogt, \emph{Sur les invariants fondamentaux des equations
  diff{\'e}rentielles lin{\'e}aires du second ordre}, Ann. Sci. d'{\'E}.N.S. 3
  s{\'e}rie \textbf{6}, 3--71 (1889).

\end{thebibliography}

\end{document}